\documentstyle[12pt]{article}
\def\a{\alpha}

\def\d{\delta}
\def\l{\lambda}

\newcommand{\bn}{\begin{equation}}
\newcommand{\ed}{\end{equation}}

\newtheorem{proposition}{Proposition}[section]
\newtheorem{theorem}{Theorem}[section]
\newtheorem{lemma}{Lemma}[section]
\newtheorem{corollary}{Corollary}[section]
\newcommand{\hsp}{\mbox{$\hspace{.3in}$}}

\begin{document}
\begin{center}
{\LARGE\bf Deformation of Yangian $Y(sl_2)$}\\
\vspace{0.5cm}
{\large\bf S.M. Khoroshkin}$^1$, {\large\bf A.A.Stolin}$^2$
{\large\bf and V.N.Tolstoy}$^3$\\
\vspace{0.2cm}
$^1$ Institute of Theoretical \& and Experimental Physics \\
117259 Moscow, Russia (e-mail:khoroshkin@vxitep.itep.ru)\\
\vspace{0.2cm}
$^2$ Department of Mathematics, Royal Institute of Technology \\
S-10044 Stockholm, Sweden (e-mail: astolin@math.kth.se)\\
\vspace{0.2cm}
$^3$ Institute of Nuclear Physics, Moscow State University \\
119899 Moscow, Russia (e-mail: tolstoy@anna19.npi.msu.su)
\end{center}
\date{}
\vspace{0.5cm}

\begin{abstract}
A quantization of a non-standard rational solution of CYBE for $sl_2$
is given explicitly. We obtain the quantization with the help of
a twisting of the usual Yangian $Y(sl_2)$.
This quantum object (deformed Yangian $Y_{\eta,\xi}(sl_2))$ is
a two-parametric deformation of the universal enveloping algebra
$U(sl_2[u])$ of the positive current algebra $sl_2[u]$.
We consider the pseudotriangular structure on $Y_{\eta,\xi}(sl_2)$,
the quantum double $DY_{\eta,\xi}(sl_2)$, its
the universal R-matrix and also the RTT-realization of
$Y_{\eta,\xi}(sl_2)$.
\end{abstract}

\setcounter{section}{-1}
\section{Introduction}

Yangian of $Y(g)$ of a simple Lie algebra $g$ was introduced by
V. Drinfeld \cite{D1,D2} as a deformation of the universal enveloping
algebra $U(g[u])$ of the current algebra $g[u]$. Recently
the Yangian symmetry $Y(sl_n)$ was shown for the following
one-dimensional N-body integrable models: the Hubbard model
\cite{UK,U}, the classical $sl_n$ Euler-Calogero-Moser model
confined in an external harmonic potential \cite{H1}, the quantum
$sl_n$ Calogero model confined in the harmonic potential \cite{H2},
and
the quantum Sutherland model \cite{HHTBP,BGHP}. Moreover the Yangian
representation theory applied to the S-matrix theory of
the $G\times G$-invariant $\sigma$-model ($G$-principal chiral model)
\cite{N}, and also to the system with non-local conserved currents
\cite{LS} and so on.

The Yangians are an useful tool not only in the physics
but in the classical representation theory of the simple
Lie algebras as well. In particular Yangians are
employed in \cite{MNO,M} for an explicit description of the center
of the universal enveloping algebra $U(g)$, where $g$ is
a simple Lie algebra of A-, B-, C-, D-series. In other words
Yangians enable one to construct new Laplace operators.
A connection between the Yangian $Y(gl_2)$  and the classical
construction of the Gelfand-Zetlin basis for the Lie algebra $gl_2$
was established in \cite{NT}.

Tensor products of finite dimensional representations of
the Yangian $Y(g)$ produce rational solutions of
the quantum Yang-Baxter equation (QYBE).
For instance, for $g=sl_2$ these solutions can be obtained by
the fusion procedure applied to the Yang solution $R(u)=1+
{\bf p}/u$, where ${\bf p}$ is the permutation of factors in
${\bf C}^2\otimes{\bf C}^2$.

Another way to obtain a rational solution of QYBE is to find
the image of the universal R-matrix for $Y(g)$ in a tensor
product of finite dimensional representations of $Y(g)$ (see
\cite{KT}). The classical r-matrix, which corresponds to the
universal R-matrix for $Y(g)$ is $r={\bf c}_2/u$, where ${\bf c}_2$ is
the canonical invariant element of $g\otimes g$, namely
${\bf c}_2=\sum_{i} I_i\otimes I^i$, where $\{I_i\}$ and
$\{I^i\}$ are dual bases for $g$ with respect to the Killing form.
Therefore, Yangians provide a "sophisticated" way to produce
rational solutions of QYBE.
However, there exist other rational solutions of the classical
Yang-Baxter equation (CYBE).  These solutions were studied in
\cite{S1}.  Every rational solution of CYBE provides a bialgebra
structure on $g[u]$.  These structures have not been quantized
yet except the case $r={\bf c}_2/u$, when the Yangian is exactly
the quantization.

We present here a quantization of the simplest non-standard
rational r-matrix for $sl_2$, namely $r={\bf c}_2/u+
h_{\alpha}\wedge e_{-\alpha}$. Here $e_{\pm\alpha},~h_{\alpha}$ is
the standard Chevalley basis for $sl_2$.  We note
that the additional term $h_{\alpha}\wedge e_{-\alpha}$ leads to
a deformation of the co-algebra structure of $Y(sl_2)$.
We perform this deformation by means of the twisting of
$U(sl_2)\otimes U(sl_2)$ by some special two-tensor of
$U(sl_2)\otimes U(sl_2)$ (appeared in \cite{GGS}).
Moreover, this two-tensor enables us to write down a quantum R-matrix
corresponding to the classical r-matrix $r={\bf c}_2/u+
h_{\alpha}\wedge e_{-\alpha}$. On the other hand, writing down
this R-matrix explicitly in the fundamental representation, we
develop RTT-formalism (see \cite{FRT}) to get another presentation of
the deformed Yangian.  We discuss also properties of the corresponding
quantum determinant, and the realization of the deformed Yangian
$Y_{\eta,\xi}(sl_2)$ in terms of generating functions
("field" realization).

The existence of the quantum determinant and the field
realization of the deformed Yangian can be used in a description
of the center of $U(g)$ and we hope to discuss this question
in future publications.

\setcounter{equation}{0}
\section{Non-Standard Quantization of $U(sl_2)$}

Let $e_{\pm\alpha},~h_{\alpha}$ be the Chevalley basis for the
universal
enveloping algebra $U(sl_2)$ of the Lie algebra $sl_2$ with the
standard
defining relations:
\bn
[e_{\alpha},e_{-\alpha}]=h_{\alpha}~,\hsp
[h_{\alpha},e_{\pm\alpha}]=\pm2e_{\pm\alpha}~.
\label{NSQ1}
\ed
Here and anywhere we put $(\alpha,\alpha)=2$. Let $U(b_-)\in U(sl_2)$
be
the universal enveloping algebra of the Borel subalgebra $b_-$ of
$sl_2$,
generated by the elements $h_{\alpha}$ and $e_{-\alpha}$.
Let us introduce the following two-tensor (a formal series) $F$ of
some
extension of $U(b_-)\otimes U(b_-)$:
\bn
F=1+\xi h_{\alpha}\otimes e_{-\alpha}+
\frac{\xi^2}{2}h_{\alpha}(h_{\alpha}+2)\otimes e_{-\alpha}^2+\ldots
=\sum_{k\geq 0}\frac{\xi^{k}}{k!}\Big(\prod_{i=0}^{k-
1}(h_{\a}+2i)\Big)
\otimes e_{-\a}^k \ ,
\label{NSQ2}
\ed
where $\xi\in {\bf C}$ is some parameter. We borrowed the element $F$
from \cite{GGS}. It is not difficult to verify that the following
series
\bn
F^{-1}=1-\xi h_{\alpha}\otimes e_{-\alpha})+
\frac{\xi^2}{2}h_{\alpha}(h_{\alpha}-2)\otimes e_{-\alpha}^2+
\ldots
=\sum_{k\geq 0}\frac{(-\xi)^{k}}{k!}\Big(\prod_{i=0}^{k-1}(h_{\a}-
2i)\Big)
\otimes e_{-\a}^k
\label{NSQ3}
\ed
is a inverse element to $F$, i.e. $FF^{-1}=F^{-1}F=1$.
The following proposition is valid.
\begin{proposition}
The element $F$ (as formal series) satisfies the following relation
\bn
F^{12}(\Delta\otimes {\rm id})F=F^{23}({\rm id}\otimes\Delta)F \ ,
\label{NSQ4}
\ed
where $\Delta$ is the usual comultiplication in $U(b_-)$, i.e.
$\Delta(a)=a\otimes 1+1\otimes a$ for any $a\in b_-$.
\label{NSP1}
\end{proposition}
{\it Proof}. By direct calculations.

As consequence of this proposition we have
\begin{corollary}
Let $\cal A$ be an arbitrary Hopf algebra containing the Hopf algebra
$U(b_-)$ and let $\bar{\cal A}_{\xi}^{(F)}$ be the algebra
${\cal A}[[\xi]]$ (i.e. the algebra ${\cal A}$ over ${\bf C}[[\xi]]$)
with co-multiplication map $\Delta^{(F)}$ given by the formula
\bn
\Delta^{(F)}(a)=F\Delta(a)F^{-1}~, \hsp
(\forall\: a\in \bar{\cal A}_{\xi}^{(F)})~,
\label{NSQ5}
\ed
then $\bar{\cal A}_{\xi}^{(F)}$ is a Hopf algebra.
\label{NSQC1}
\end{corollary}
{\it Proof}. Coassociativity of $\Delta^{(F)}$ follows from the
formula (\ref{NSQ4}). Existence of the antipod is proved in \cite{D3}.

Now we introduce the following notations for some elements of
$\bar{U}_{\xi}(b_-)^{(F)}$:
\bn
T_{\alpha}:=1-\xi e_{-\alpha}~, \hsp
T_{\alpha}^{-1}:=(1-\xi e_{-\alpha})^{-1}~.
\label{NSQ6}
\ed
\begin{proposition}
The elements $h_{\alpha}$, $T_{\alpha}^{\pm 1}$ satisfy the following
relations:
\bn
T_{\alpha}T_{\alpha}^{-1}=T_{\alpha}^{-1}T_{\alpha}~,\hsp
[h_{\alpha},T_{\alpha}]=2(1-T_{\alpha})~, \hsp
[h_{\alpha},T_{\alpha}^{-1}]=2(T_{\alpha}^{-1}-T_{\alpha}^{-2})~,
\label{NSQ7}
\ed
\bn
\Delta^{(F)}(h_{\alpha})=h_{\alpha}\otimes T_{\alpha}^{-1}+
1\otimes h_{\alpha}~, \hsp
\Delta^{(F)}(T_{\alpha})=T_{\alpha}\otimes T_{\alpha}~, \hsp
\Delta^{(F)}(T_{\alpha}^{-1})=T_{\alpha}^{-1}\otimes T_{\alpha}^{-1},
\label{NSQ8}
\ed
\bn
S(h_{\alpha})=-h_{\alpha}T_{\a}~, \hsp S(T_{\alpha})=T_{\alpha}^{-
1}~,
\hsp S(T_{\alpha}^{-1})=T_{\alpha}~,
\label{NSQ9}
\ed
\bn
\varepsilon(h_{\alpha})=0~, \hsp
\varepsilon(T_{\alpha})=\varepsilon(T_{\alpha}^{-1})=1~.
\label{NSQ10}
\ed
The algebra $U_{\xi}(b_-)_{(F)}$ generated by $h_{\alpha}$,
$T_{\alpha}^{\pm 1}$ is a Hopf subalgebra of
$\bar{U}_{\xi}(b_-)^{(F)}$.
\label{NSP2}
\end{proposition}
{\it Proof}. By direct calculations with (\ref{NSQ1}), (\ref{NSQ4})
and
(\ref{NSQ5}).

The Hopf algebra $U_{\xi}(b_-)^{(F)}$ is triangular with the universal
R-matrix
\[
R=F^{21}F^{-1}=1+\xi e_{-\alpha}\wedge h_{\alpha}+\ldots =
\]
\bn
=\sum_{k,m\geq 0}\frac{(-1)^m\xi^{k+m}}{k!m!}
\Big(\prod_{i=0}^{m-1}(h_{\a}+2k-2i)\Big)e_{-\a}^k\otimes
\Big(\prod_{i=0}^{k-1}(h_{\a}+2j)\Big)e_{-\a}^m =
\label{NSQ11}
\ed
\[
=\sum_{k,m\geq 0}\frac{(-1)^m\xi^{k+m}}{k!m!}
\Big(\prod_{i=0}^{m-1}(h_{\a}+2k-2i)\Big)(1-T_{\a})^k\otimes
\Big(\prod_{i=0}^{k-1}(h_{\a}+2j)\Big)(1-T_{\a})^m =
\]
\bn
R^{21}=R^{-1}~.
\label{NSQ12}
\ed
This algebra is a quantization of the bialgebra Lie $b_-$ defined
by the classical r-matrix $r=e_{-\a}\wedge h_{\a}$.

\noindent
According to Corollary \ref{NSQC1} we can extend the twisting by $F$
to $\bar{U}_{\xi}(sl_2)^{(F)}$. Then we have
\begin{proposition}
Let $U_{\xi}(sl_2)^{(F)}$ be an algebra generated the elements
$h_{\a}$, $e_{\a}$  and $T_{\a}^{\pm1}$ with the defining relations
\bn
T_{\alpha}T_{\alpha}^{-1}=T_{\alpha}^{-1}T_{\alpha}~,\hsp
[h_{\alpha},T_{\alpha}]=2(1-T_{\alpha})~, \hsp
[h_{\alpha},T_{\alpha}^{-1}]=2(T_{\alpha}^{-1}-T_{\alpha}^{-2})~,
\label{NSQ13}
\ed
\bn
[h_{\a},e_{\a}]=2e_{\a}~,\hsp [T_{\a}, e_{\a}]=2\xi h_{\a}~,\hsp
[T_{\a}^{-1}, e_{\a}]=-2\xi T_{\a}^{-1}h_{\a}T_{\a}^{-1}~,
\label{NSQ14}
\ed
\bn
\Delta^{(F)}(h_{\alpha})=h_{\alpha}\otimes T_{\alpha}^{-1}+
1\otimes h_{\alpha}~, \hsp
\Delta^{(F)}(T_{\alpha})=T_{\alpha}\otimes T_{\alpha}~, \hsp
\Delta^{(F)}(T_{\alpha}^{-1})=T_{\alpha}^{-1}\otimes T_{\alpha}^{-1},
\label{NSQ15}
\ed
\bn
\Delta^{(F)}(e_{\a})=e_{\a}\otimes T_{\a}^{-1}+1\otimes e_{\a}-
\xi h_{\a}\otimes T_{\a}^{-1}h_{\a}-
\frac{\xi}{2}h_{\a}(h_{\a}-2)\otimes T_{\a}^{-1}-
\frac{\xi}{2}h_{\a}(h_{\a}+2)\otimes T_{\a}^{-2}~,
\label{NSQ16}
\ed
\bn
S(h_{\alpha})=-h_{\alpha}T_{\a}~, \hsp
S(T_{\alpha})=T_{\alpha}^{-1}~, \hsp
S(T_{\alpha}^{-1})=T_{\alpha}~,
\label{NSQ17}
\ed
\bn
S(e_{\a})=-e_{\a}T_{\a}-\frac{\xi}{2}h_{\a}(h_{\a}+2)T_{\a}(T_{\a}-
2)~,
\label{NSQ18}
\ed
\bn
\varepsilon(h_{\alpha})=\varepsilon(e_{\alpha})=0~, \hsp
\varepsilon(T_{\alpha})=\varepsilon(T_{\alpha}^{-1})=1~.
\label{NSQ19}
\ed
Then $U_{\xi}(sl_2)^{(F)}$ is a Hopf subalgebra of
$\bar{U}_{\xi}(sl_2)^{(F)}$ and it is a triangular deformation
of $U(sl_2)$ in the direction of the classical r-matrix
$r=e_{-\a}\wedge h_{\a}$.
\label{NSP3}
\end{proposition}
{\it Proof}. By direct calculations with (\ref{NSQ1}) and
(\ref{NSQ5}).

\noindent
{\it Remark 1}. In every finite dimensional representation of
$sl_2$ the element $T_{\a}=1-2\xi e_{-\a}$ is always
invertible since $e_{-\a}$ is nilpotent. Therefore, the theory of
finite dimensional representations of $U_{\xi}(sl_2)^{(F)}$ is the
same as the theory for $sl_2$.

\noindent
{\it Remark 2}. Similar computations for $sl_2$ with another
twisting element $\tilde{F}$ were carried out in \cite{O}.
However, using Theorem 2 from \cite{D4} one can prove the following
result.
\begin{theorem}
There exists an invertible element $T\in U(gl_2)[[\xi]]$ such that
$\tilde{F}=(T~\otimes~T)\cdot F^{-1}\cdot$$\Delta(T^{-1})$,
$\varepsilon(T)=1$ and ${\rho}_1 (T)=1\in GL_2$, where  ${\rho}_1$ is
the two-dimensional representation of $gl_2$.
In other words Hopf algebras obtained from $U(sl_2)[[\xi]]$
by twisting by  $\tilde{F}^{-1}$ and ${F}$ are isomorphic
as Hopf algebras and the isomorphism is given by conjugation
by $T$.
\label{NSQT1}
\end{theorem}
{\it Proof}. The matrix $({\rho}_1\otimes
{\rho}_1)((\tilde{F}^{21})^{-1}\tilde{F})$
was computed in \cite{O}. We can calculate the matrix
$({\rho}_1\otimes{\rho}_1) ({F}^{21}F^{-1})$ and it turns out
that we obtain the same matrix (see Lemma 5.1 further).
Then Theorem 2 from \cite{D4} implies all the statements.

\setcounter{equation}{0}
\section{Rational solutions of CYBE}

The theory of rational solutions of CYBE over any simple Lie algebra
$g$ was developed in \cite{S1} and we would like to remind several
basic facts of the theory.

Let $P(u,v)=\frac{{\bf c}_2}{u-v}+r(u,v)$ be a function from $
{\bf C}^2$ to $g\otimes g$, where $r(u,v)$ is a polynomial in $u,v$.
If $P(u,v)$ satisfies CYBE, we say that $P(u,v)$ is
a rational solution of CYBE.  It is possible to show that every
rational
solution can be brought by means of a gauge transformation to the form
\bn
P(u,v)=\frac{{\bf c}_2}{u-v}+r_{00}+r_{10}u+r_{01}v+r_{11}uv~,
\label{RS1}
\ed
where $r_{00},~r_{10},~r_{01},~r_{11}~\in g\otimes g$.
Every rational solution induces a bialgebra structure on the
current algebra $g[t]$.  It is possible to prove the following
theorems.
\begin{theorem}
Let $D(g[t])$ be the classical double corresponding to a rational
solution of CYBE.  Then $D(g[t])$ and $g((t^{-1}))$ are isomorphic
as Lie algebras with inner product which takes the following form on
$g((t^{-1}))$:
\bn
(a(t),b(t))=Res_{t=0}\langle a(t),b(t)\rangle~.
\label{RS2}
\ed
where $a(t),~b(t),~\in g((t^{-1}))$. It means that $g((t^{-1}))$
can be represented as a Manin triple $g((t^{-1}))=g[t]\oplus {\cal W}$
where ${\cal W}$ is a Lagrangian subalgebra with respect to the inner
product
(\ref{RS2}).
\label{RST1}
\end{theorem}
\begin{theorem}
Let $P(u,v)=\frac{{\bf c}_2}{u-v}+r_{00}$, where $r_{00}\in g\otimes
g$.
Then the following conditions are equivalent:

\noindent
{\it(i)} the function $P(u,v)$ is a rational r-matrix;

\noindent
{\it (ii)} the element $r_{0,0}$ satisfies CYBE;

\noindent
{\it (iii)} the Lagrangian subalgebra ${\cal W}$ is contained in
$g[[t^{-1}]]$ (this is equivalent to $t^{-2}g[[t^{-1}]]\subset
{\cal W}\subset g[[t^{-1}]]$).
\label{RST2}
\end{theorem}
Now it is possible to deduce that rational solutions satisfying
the Theorem (\ref{RST2}) are in a 1-1 correspondence with the
following
combinatorial data:

\noindent
{\it 1)} subalgebra ${\cal L}$ of $g$;

\noindent
{\it 2)} non-degenerate 2-cocycle $B$ on ${\cal L}$.

In case of $g=sl_n$ all the rational solutions can be described in a
similar way. Let
\bn
\hbox{$d_k=diag(~$\vtop to 28pt{\baselineskip 6pt  
\hbox to 39pt{\hfill$1,\ldots,1$\hfill}
\hbox to 39pt{\upbracefill}
\vskip 3pt
\hbox to 39pt{\hfill$k$\hfill}\vfill}
$t,\ldots,t)\in GL(n,{\bf C}((t^{-1})))~.$}
\label{RS3}
\ed
Then every rational solution of CYBE defines some Lagrangian
subalgebra
${\cal W}$ contained in $d_{k}^{-1}\cdot sl(n,{\bf C}[[t^{-1}]])\cdot
d_k$
for some $k$. The corresponding combinatorial data are:

\noindent
(1) subalgebra ${\cal L}\subset sl(n,{\bf C})$ such that
${\cal L}+{\cal P}_k=sl(n,{\bf C})$, where ${\cal P}_k$ is
the maximal parabolic subalgebra of $sl(n,{\bf C})$ not containing
the root vector $e_{\a_k}$ of the simple root $\a_{k}$;

\noindent
(2) 2-cocycle $B$ on ${\cal L}$ which is nondegenerate on
${\cal L}\cap {\cal P}_k$.

In case of $sl_2$ one has just two non-standard rational r-matrices
(up to gauge equivalence):
\bn
P_1(u,v)=\frac{{\bf c}_2}{u-v}+h_{\a}\wedge e_{-\a}~,
\label{RS4}
\ed
and
\bn
P_2(u,v)=\frac{{\bf c}_2}{u-v}+h_{\a}\otimes e_{-\a}u-e_{-\a}\otimes
h_{\a}v~.
\label{RS5}
\ed
The corresponding Lagrangian subalgebras are
\bn
{\cal W}_1=t^{-2}(sl_2[[t^{-1}]])\oplus{\bf C}\cdot(e_{\a}t^{-1}-
h_{\a})
\oplus {\bf C}(e_{-\a}t^{-1})\oplus {\bf C}(h_{\a}t^{-1}+2e_{-\a})~,
\label{RS6}
\ed
and
\[
{\cal W}_2=t^{-3}(sl_2[[t^{-1})\oplus {\bf C}(e_{\a}t^{-1}-h_{\a}t)
\oplus {\bf C}(e_{-\a}t^{-1})\oplus
\]
\bn
\oplus{\bf C}(h_{\a}t^{-1})
\oplus {\bf C}(e_{\a}t^{-2})\oplus {\bf C}(e_{-\a}t^{-2})
\oplus {\bf C}(h_{\a}t^{-2}+2e_{-\a})~.
\label{RS7}
\ed
Now we are going to quantize the rational r-matrix (\ref{RS4}).

\setcounter{equation}{0}
\section{Twisting of Yangian $Y(sl_2)$}

One can show that the Yangian $Y(sl_2)$ (as a Hopf algebra) can be
defined by Chevalley generators $h_{\a},~e_{\pm\a},~e_{\d-\a}$
with the defining relations \cite{KT}:
\bn
[e_{\alpha},e_{-\alpha}]=h_{\alpha}~,\hsp
[h_{\alpha},e_{\pm\alpha}]=\pm2e_{\pm\alpha}~,
\label{TY1}
\ed
\bn
[h_{\alpha},e_{\d-\alpha}]=e_{\d-\alpha}~,\hsp
[e_{-\alpha},e_{\d-\alpha}]=\eta e_{-\alpha}~,
\label{TY2}
\ed
\bn
[e_{\a},[e_{\a},[e_{\a},e_{\d-\alpha}]]]=6\eta e_{\a}^{2}~,\hsp
[[[e_{\a},e_{\d-\a}],e_{\d-\a}],e_{\d-\a}]=6\eta e_{\d-\a}^{2}~.
\label{TY3}
\ed
\bn
\Delta(h_{\a})=h_{\a}\otimes 1 + 1\otimes h_{\a}~, \hsp
\Delta(e_{\pm\a})=e_{\pm\a}\otimes 1 + 1\otimes e_{\pm\a}~,
\label{TY4}
\ed
\bn
\Delta(e_{\d-\a})=e_{\d-\a}\otimes 1 + 1\otimes e_{\d-\a}+
\eta e_{\a}\otimes h_{\a}~,
\label{TY5}
\ed
\bn
S(h_{\a})=-h_{\a}~, \hsp S(e_{\pm\a})=-e_{\pm\a}~,
\hsp S(e_{\d-\a})=-e_{\d-\a}-\eta e_{-\a}h_{\a}~,
\label{TY6}
\ed
\bn
\varepsilon (h_{\a})=\varepsilon (e_{\pm\a})=\varepsilon (e_{\d-
\a})=0~,
\hsp \varepsilon (1)~=1~.
\label{TY7}
\ed
\noindent
Here we explicitly introduce the Yangian deformation parameter
$\eta \in {\bf C}\backslash \{0\}${\footnote{All Yangians with
different
$\eta$ are isomorphic one to another and therefore one usually suppose
that $\eta=1$.} and therefore we shall use further the notation
$Y_{\eta}(sl_2)$ for the Yangian.

Since $Y_{\eta}(sl_2)$ contains $U(sl_2)$ as a Hopf subalgebra
the Corollary \ref{NSQC1} implies that the algebra
${\bar Y}_{\eta,\xi}^{(F)}$ isomorphic to $Y_{\eta}[[\xi]]$
with the comultiplication (\ref{NSQ5}) is a Hopf algebra.
\begin{proposition}
The elements $h_{\a}~, e_{\a}~, T_{\a}^{\pm 1}$ (see (\ref{NSQ6})) and
$e_{\d-\a}$ satisfy the relations:
\bn
T_{\alpha}T_{\alpha}^{-1}=T_{\alpha}^{-1}T_{\alpha}~,\hsp
[h_{\alpha},T_{\alpha}]=2(1-T_{\alpha})~, \hsp
[h_{\alpha},T_{\alpha}^{-1}]=2(T_{\alpha}^{-1}-T_{\alpha}^{-2})~,
\label{TY8}
\ed
\bn
[h_{\a},e_{\a}]=2e_{\a}~,\hsp [T_{\a}, e_{\a}]=2\xi h_{\a}~,\hsp
[T_{\a}^{-1}, e_{\a}]=-2\xi T_{\a}^{-1}h_{\a}T_{\a}^{-1}~,
\label{TY9}
\ed
\bn
[T_{\a}, e_{\d-\a}]=-\frac{\eta}{2\xi}\left(T_{\a}^{2}-
2T_{\a}+1\right)~,
\hsp [T_{\a}^{-1}, e_{\d-\a}]=-\frac{\eta}{2\xi}
\left(T_{\a}^{-2}-2T_{\a}^{-1}+1\right)~,
\label{TY10}
\ed
\bn
[e_{\a},[e_{\a},[e_{\a},e_{\d-\alpha}]]]=6\eta e_{\a}^{2}~,\hsp
[[[e_{\a},e_{\d-\a}],e_{\d-\a}],e_{\d-\a}]=6\eta e_{\d-\a}^{2}~.
\label{TY11}
\ed
\bn
\Delta^{(F)}(h_{\alpha})=h_{\alpha}\otimes T_{\alpha}^{-1}+
1\otimes h_{\alpha}~, \hsp
\Delta^{(F)}(T_{\alpha})=T_{\alpha}\otimes T_{\alpha}~, \hsp
\Delta^{(F)}(T_{\alpha}^{-1})=T_{\alpha}^{-1}\otimes T_{\alpha}^{-1},
\label{TY12}
\ed
\bn
\Delta^{(F)}(e_{\a})=e_{\a}\otimes T_{\a}^{-1}+1\otimes e_{\a}-
\xi h_{\a}\otimes T_{\a}^{-1}h_{\a}-
\frac{\xi}{2}h_{\a}(h_{\a}-2)\otimes T_{\a}^{-1}-
\frac{\xi}{2}h_{\a}(h_{\a}+2)\otimes T_{\a}^{-2}~,
\label{TY13}
\ed
\bn
\Delta^{(F)}(e_{\d-\a})=e_{\a}\otimes T_{\a}+1\otimes e_{\d-\a}+
\xi h_{\a}\otimes T_{\a}^{-1}+
\frac{\eta}{2\xi}h_{\a}\otimes (1-T_{\a})~,
\label{TY14}
\ed
\bn
S(h_{\alpha})=-h_{\alpha}T_{\a}~, \hsp S(T_{\alpha})=T_{\alpha}^{-1}~,
\hsp S(T_{\alpha}^{-1})=T_{\alpha}~,
\label{TY15}
\ed
\bn
S(e_{\a})=-e_{\a}T_{\a}-\frac{\xi}{2}h_{\a}(h_{\a}+2)T_{\a}(T_{\a}-
1)~,
\label{TY16}
\ed
\bn
S(e_{\d-\a})=-e_{\d-\a}T_{\a}^{-1}- \frac{\xi}{\eta}h_{\a}T_{\a}+
\frac{\eta}{2\xi}h_{\a}T_{\a}^{-1} -\frac{\eta}{2\xi}h_{\a}~,
\label{TY17}
\ed
\bn
\varepsilon(h_{\alpha})=\varepsilon(e_{\alpha})=
\varepsilon(e_{\d-\a})=0~, \hsp
\varepsilon(T_{\alpha})=\varepsilon(T_{\alpha}^{-1})=1~.
\label{TY18}
\ed
The algebra $Y_{\eta,\xi}(sl_2)$ generated the elements
$h_{\a}$, $e_{\a}$, $T_{\a}^{\pm 1}$, $e_{\d-\a}$
is a Hopf subalgebra of $\bar{Y}_{\eta,\xi}(sl_2)^{(F)}$.
\label{P2}
\end{proposition}
{\it Proof}. By direct calculations.

One can see that the Hopf algebra $Y_{\eta,\xi}sl_2)^{(F)}$ is a
quantization of the Lie bialgebra $sl_2[u]$ corresponding to
the rational solution (\ref{RS4}).

\setcounter{equation}{0}
\section{Pseudotriangular structure on $Y_{\eta,\xi}(sl_2)$.
Deformed Yangian double $DY_{\eta}(sl_2)$}

We want to recall that $Y_{\eta}(sl_2)$ is a pseudotriangular
Hopf algebra \cite{D1}. It means that there is a family of
automorphisms
${\cal T}_\l:Y_{\eta}(sl_2)\mapsto Y_{\eta}(sl_2)$ and an element
$R(\l)=1+\sum^\infty_{k=1}R_k\l^{-1}$, where
$R_k\in Y_{\eta}(sl_2)\otimes Y_{\eta}(sl_2)$, such that
\bn
({\cal T}_\l\otimes {\cal T}_\mu)R(u)=R(u+\l-\mu)~, \hsp
({\cal T}_\l\otimes {\rm id})\Delta(a)'=R(\l)(({\cal T}_\l
\otimes {\rm id})\Delta(a))R(\l)^{-1}~,
\label{PS1}
\ed
\bn
(\Delta\otimes {\rm id})R(\l)=R^{13}(\l)R^{23}(\l)~, \hsp
R^{12}(\l)R^{21}(-\l)=1\otimes 1~,
\label{PS2}
\ed
\bn
R^{12}(\l_1-\l_2)R^{13}(\l_1-\l_3)R^{23}(\l_2-\l_3)=
R^{23}(\l_2-\l_3)R^{13}(\l_1-\l_3)R^{12}(\l_1-\l_2)~.
\label{PS3}
\ed
As was proved in \cite{S2}, ${\bar Y}_{\eta,\xi}(sl_2)^{(F)}$ is also
a pseudotriangular Hopf algebra with respect to $\Delta^{(F)}=
F\Delta F^{-1}$, $R(\l)^{(F)}=F^{21}R(\l)F^{-1}$ and
the same ${\cal T}_{\l}$.

Since $Y_{\eta}(sl_2)\subset Y_{\eta,\xi}(sl_2)
\subset{\bar Y}_{\eta,\xi}(sl_2){(F)}$ (as associative algebras)
and since ${\cal T}_\l$ acts
identically on $U(sl_2)[[\xi]]$ and since $Y_{\eta,\xi}(sl_2)$ differs
from $Y_{\eta}(sl_2)$ by elements from $U(sl_2)[[\xi]]$ therefore
the automorphisms ${\cal T}_\l$ extended to
${\bar Y}_{\eta,\xi}(sl_2)^{(F)}$ preserve $Y_{\eta,\xi}(sl_2)$.
Thus we have the following result.
\begin{proposition} The Hopf algebra $Y_{\eta,\xi}(sl_2)$ is
pseudotriangular with $R(\l)^{(F)}:=F^{21}R(\l)F^{-1}$. In
particular, $R(\l)^{(F)}$ is a rational solution of QYBE.
\label{PSP1}
\end{proposition}

Let us recall that the Yangian double $DY_{\eta}(sl_2)$ (see
\cite{KT})
is a quasitriangular Hopf algebra with an universal R-matrix ${\cal
R}$
which lies in some extension of $DY_{\eta}(sl_2)\otimes
DY_{\eta}(sl_2)$.
Since $U(sl_2)\subset DY_{\eta}(sl_2)$ we can twist the Yangian double
$DY_{\eta}(sl_2)$ by $F$. Using formal algebraic arguments similar to
that of \cite{S2} we get as result the following proposition.
\begin{proposition} The deformed Yangian double
$DY_{\eta,\xi}(sl_2)^{(F)}$ is a quasitriangular Hopf algebra
with the universal R-matrix ${\cal R}^{(F)}=F^{21}{\cal R}F^{-1}$.
\label{PSP2}
\end{proposition}

In what follows we need the realization of the Yangian double
$DY_{\eta}(sl_2)$ given in \cite{KT}. In this realization the Yangian
double $DY_{\eta}(sl_2)$ is generated by the elements $h_{k\d}$,
$e_{k\d\pm\a}$, ($k\in {\bf Z}$)\footnote{These notations of
the generators are connected with the notations in \cite{KT} as
follows:
$h_{k\d}:=h_k$, $e_{k\d+\a}:=e_k$, $e_{k\d-\a}:=f_k$ ($k\in {\bf
Z}$).},
which are composed into generating functions
$h^+(u)=1+\sum_{k\ge0}h_{k\d}u^{-k-1}$, $(h_0:=h_{\a})$,
$e_{\pm\a}^+(u)=\sum_{k\ge0}e_{k\d\pm\a}u^{-k-1}$,
and $h^-(u)=1-\sum_{k<0}h_{k\d}u^{-k-1}$,
$e_{\pm\a}^-(u)=-\sum_{k<0}e_{k\d\pm\a}u^{-k-1}$,
which satisfy the following relations:
\bn
[h^{\pm}(u),h^{\pm}(v)]=0~,\hsp [h^+(u),h^-(v)]=0~,
\label{PS4}
\ed
\bn
[e_{\a}^{\pm}(u),e_{-\a}^{\pm}(v)]=
-\eta\frac{h^{\pm}(u)-h^{\pm}(v)}{u-v}~,
\label{PS5}
\ed
\bn
[e_{\a}^{\pm}(u),e_{-\a}^{\mp}(v)]=
-\eta\frac{h^{\mp}(u)-h^{\pm}(v)}{u-v}~,
\label{PS6}
\ed
\bn
[h^{\pm}(u),e_{\pm\a}^{\pm}(v)]=
-\eta\frac{\{h^{\pm}(u),(e_{\pm\a}^{\pm}(u)-
e_{\pm\a}^{\pm}(v))\}}{u-v}~,
\label{PS7}
\ed
\bn
[h^{\pm}(u),e_{\pm\a}^{\mp}(v)]=
-\eta\frac{\{h^{\pm}(u),(e_{\pm\a}^{\pm}(u)-e_{\pm\a}^{\mp}(v))\}}{u-
v}~,
\label{PS8}
\ed
\bn
[e_{\pm\a}^{\pm}(u),e_{\pm\a}^{\pm}(v)]=
\mp\eta\frac{(e_{\pm\a}^{\pm}(u)-e_{\pm\a}^{\pm}(v))^2}{u-v}~,
\label{PS9}
\ed
\bn
[e_{\pm\a}^+(u),e_{\pm\a}^-(v)]=\mp\eta\frac{(e^+(u)-e^-(v))^2}{u-v}~,
\label{PS10}
\ed
where $\{a,b\}=ab+ba$.  It turns out that $Y_{\eta}(sl_2)$ is
generated by $h^+(u)$, $e_{\pm\a}^+(u)$, while the dual to
$Y_{\eta}(sl_2)$ algebra $Y_{\eta}^{\circ}(sl_2)$ is generated by
$h^-(u)$, $e_{\pm\a}^-(u)$. The universal R-matrix ${\cal R}$ was
found
in \cite{KT} can be factorized as follows:
\bn
{\cal R}={\cal R}_+{\cal R}_0{\cal R}_-~,
\label{PS11}
\ed
where
\bn
{\cal R}_+=\prod^\rightarrow_{k\ge0}\exp\Bigl(-e_{-(k+1)\d-\a}\otimes
e_{k\d+\a}\Bigr)~, \hsp
{\cal R}_-=\prod^\leftarrow_{k\ge0}\exp\Bigl(-e_{-(k+1)\d+\a}\otimes
e_{k\d-\a}\Bigr)~,
\label{PS12}
\ed
\bn
{\cal R}_0=\prod_{n\ge0}\exp {\rm Res}_{u=v}\Bigl(\ln\,h^-
(v+2n+1)\otimes
\frac{{\rm d}}{{\rm d}u}\ln\,h^+(u)\Bigr)~.
\label{PS13}
\ed
Here ${\rm Res}_{u=v}(f(u)\otimes g(v))=\sum_k f_k\otimes g_{-k-1}$
if $f(u)=\sum f_ku^{-k-1}$, $g(v)=\sum g_kv^{-k-1}$.
\begin{corollary}
The element ${\cal R}^{(F)}=F{\cal R}F^{-1}$ satisfies QYBE,
where $F$ is the same as in (\ref{NSQ2}) and ${\cal R}$ is defined by
(\ref{PS11})-(\ref{PS13}).
\end{corollary}

\setcounter{equation}{0}
\section{RTT-realization of the deformed Yangian $Y_{\eta,\xi}(sl_2)$}

In the section we develop so called RTT-formalism (see \cite{FRT}),
i.e. we obtain the RTT-realization or in other words the realization
in terms of the L-operator.

Let $\rho^{(1)}$ be the two-dimensional representation of $sl_2$
in ${\bf C}^2$ with the basis $|1\rangle$ and $|-1\rangle$.  It is
well-known
that $\rho^{(1)}$ is extended to a representation  $\rho^{(1)}_u$
($u\in {\bf C}$) of the Yangian $Y_{\eta}(sl_2)$ by means of
\[
\rho^{(1)}_{u}(h^\pm(w))\,|1\rangle=\frac{\eta}{w-u}\,|1\rangle~, \hsp
\rho^{(1)}_{u}(h^\pm(w))\,|-1\rangle=\frac{-\eta}{w-u}\,|-1\rangle~.
\]
\bn
\rho^{(1)}_{u}(e_{\a}^\pm(w))\,|1\rangle=0~,\hsp
\rho^{(1)}_{u}(e_{\a}^\pm(w))\,|-1\rangle=\frac{\eta}{w-
u}\,|1\rangle~,
\label{RTT1}
\ed
\[
\rho^{(1)}_{u}(e_{-\a}^\pm(w))\,|1\rangle=\frac{\eta}{w-u}\,|-
1\rangle~,
\hsp \rho^{(1)}_{u}(e_{-\a}^\pm(w))\,|-1\rangle=0~,
\]
With the help of these formulas we find
\bn
(\rho^{(1)}_{u}\otimes\rho^{(1)}_{v})({\cal R})=
\varphi(u-v)\Bigr(1+\frac{{\bf p}_{12}}{u-v}\Bigl)=\varphi(u-v) R~(u-
v),
\label{RTT2}
\ed
where $\varphi$ is a scalar function and ${\bf p}_{12}$ interchanges
factors in ${\bf C}^2\otimes {\bf C}^2$.

Let $L(u)=(\rho^{(1)}_{u}\otimes{\rm id})({\cal R})$ then
QYBE for $R(u-v)$ implies that
\bn
R(u-v)\stackrel{1}{L}(u)\stackrel{2}{L}(v)=
\stackrel{2}{L}(v)\stackrel{1}{L}(u)R(u-v)~,
\label{RTT3}
\ed
where $\stackrel{1}{L}(u)=L(u)\otimes{\rm id}$,
$\stackrel{2}{L}(v)={\rm id}\otimes L(v)$.
The matrix $L(u)$ is a generating function for $Y_{\eta}(gl_2)$
and $Y(sl_2)\cong Y_{\eta}(gl_2)/({\rm qdet}L(u)-1)$.
More exactly we can formulate the following result.
\begin{proposition}
Let $L(u)$ be a $2\times 2$-matrix with noncommuting entries, such
that

\noindent
{\rm (i)} $R(u-v)\stackrel{1}{L}(u)\stackrel{2}{L}(v)=
\stackrel{2}{L}(v)\stackrel{1}{L}(u)R(u-v)$~,

\noindent
{\rm (ii)} $L(u)=1+\frac{L_{(0)}}{u}+\frac{L_{(1)}}{u^2}+\cdots
\frac{L_{(k)}}{u^k}+\cdots $~~,

\noindent
{\rm (iii)} ${\rm qdet}~L(u)=e_{11}(u)e_{22}(u-1)-e_{21}(u) e_{12}(u-
1)=
e_{22}(u)e_{11}(u-1)-e_{12}(u)e_{21}(u-1)=1$.

\noindent
Then the matrix coefficients of $L(u)$ generate a Hopf algebra
isomorphic to $Y_{\eta}(sl_2)$. The comultiplication $\Delta$
and the antipode $S$ are given by the formulas
\bn
\Delta(e_{ij}(u))=\sum_k e_{ik}(u)\otimes e_{kj}(u)~.
\label{RTT4}
\ed
\bn
S(L(u))=L^{-1}(u)~.
\label{RTT5}
\ed
\label{RTTP1}
\end{proposition}
The deformed Yangian $Y_{\eta,\xi}(sl_2)$ admits a similar
representation.
We start with the lemma.
\begin{lemma}
In the representation $\rho^{(1)}_u\otimes \rho^{(1)}_v$ the universal
R-matrix ${\cal R}^{(F)}=F^{21}{\cal R}F^{-1}$ has the form
\[
R_{\eta,\xi}(u-v):=(\rho^{(1)}_u\otimes\rho^{(1)}_v)({\cal R}^{(F)})=
\]
\[
=\Bigl(1+\xi\rho^{(1)}_u(e_{-\a})\otimes
\rho^{(1)}_v(h_{\a})\Bigr)\Bigl(1-\eta \frac{{\bf p}_{12}}{u-v}\Bigr)
\Bigl(1-\xi \rho^{(1)}_u(h_{\a})\otimes \rho^{(1)}_v(e_{-\a})\Bigr)=
\]
\bn
=\left(\begin{array}{cccc}
1-\frac{\eta}{u-v}& 0 & 0 & 0\\
-\xi & 1 & -\frac{\eta}{u-v} & 0\\
\xi &-\frac{\eta}{u-v} & 1 & 0\\
\xi^2 & -\xi & \xi & 1-\frac{\eta}{u-v}\end{array}\right)
\label{RTT6}
\ed
\label{RTTL1}
\end{lemma}
{\it Proof}. By direct calculation.

Let us consider an algebra $A$ of matrix elements of $L(u)$
satisfying the relation
\bn
R_{\eta,\xi}(u-v)\stackrel{1}{L}(u)\stackrel{2}{L}(v)=
\stackrel{2}{L}(v)\stackrel{1}{L}(u)R_{\eta,\xi}(u-v)~.
\label{RTT7}
\ed
It follows that $L(u)=(\rho^{(1)}_u\otimes {\rm id}){\cal R}^{(F)}$
satisfies (\ref{RTT7}). Algebra $A$ together with the comultiplication
(\ref{RTT4}) and antipod (\ref{RTT5}) constitutes a Hopf algebra.
The following lemma takes place.
\begin{lemma} {\rm (i)} The matrix $R_{\eta,\xi}(\eta)$ is the
projector
onto the one-dimensional subspace ${\bf C}
(|1\rangle\otimes |-1\rangle -
|-1\rangle \otimes |1\rangle-\xi |-1\rangle\otimes |-1\rangle)$
up to a scalar factor.

\noindent
{\rm (ii)} The following relation hold
\bn
R_{\eta,\xi}(\eta)\stackrel{1}{L}(u)\stackrel{2}{L}(u-\eta)=
\stackrel{2}{L}(u-\eta)\stackrel{1}{L}(u)R_{\eta,\xi}(\eta)=
\Bigl({\rm qdet}_{\eta,\xi}\,L(u)\Bigr)R_{\eta,\xi}(\eta)~,
\label{RTT8}
\ed
\[
{\rm qdet}_{\eta,\xi}\, L(u)=e_{11}(u)e_{22}(u-\eta)-
e_{21}(u)e_{12}(u-\eta)-\xi e_{11}(u)e_{12}(u-\eta)=
\]
\bn
=e_{22}(u)e_{11}(u-\eta)-e_{12}(u)e_{21}(u-\eta)+
\xi e_{12}(u)e_{11}(u-\eta)~,
\label{RTT9}
\ed
\bn
\Delta^{(F)}({\rm qdet}_{\eta,\xi}\, L(u))=
{\rm qdet}_{\eta,\xi}\, L(u)\otimes {\rm qdet}_{\eta,\xi}\, L(u)
\label{RTT10}
\ed
where the quantum determinant ${\rm qdet}_{\eta,\xi} L(u)$ is
an element of the Hopf algebra $A$.
\label{RTTL2}
\end{lemma}
{\it Proof}. The part (i) is verified by direct calculations.
The proof of the second part is standard (see \cite{MNO}).
\begin{lemma}
The quantum determinant ${\rm qdet}_{\eta,\xi}\, L(u)$ is
a central element of the algebra $A$.
\label{RTTL3}
\end{lemma}
{\it Proof}. The formula for the quantum determinant
${\rm qdet}_{\eta,\xi}$ and the quantum Yang-Baxter equation for
$R_{\eta,\xi}(u)$ provide the following equality:
$$
R_{\eta,\xi}^{23}(\eta)R_{\eta,\xi}^{12}(u)R_{\eta,\xi}^{13}(u+\eta)
R_{\eta,\xi}^{23}(\eta)\stackrel{1}{L}(v)\,({\rm qdet}_{\eta,\xi}\,
L(u))
R_{\eta,\xi}^{23}(\eta)=
$$
$$
=({\rm qdet}_{\eta,\xi}\,
L(u))R_{\eta,\xi}^{23}(\eta)\stackrel{1}{L}(v)
R_{\eta,\xi}^{23}(\eta)R_{\eta,\xi}^{12}(u)R_{\eta,\xi}^{13}(u+\eta)
R_{\eta,\xi}^{23}(\eta)~,
$$
where we use the standard notations:
$\stackrel{1}{L}(u):=L(u)\otimes 1\otimes 1$, and
$R_{\eta,\xi}^{12}=\sum_ia_i\otimes b_i\otimes 1$ and so on if
$R_{\eta,\xi}=\sum_ia_i\otimes b_i$.
Direct calculations show that
$$
R_{\eta,\xi}^{23}(\eta)R_{\eta,\xi}^{12}(u)R_{\eta,\xi}^{13}(u+\eta)
R_{\eta,\xi}^{23}(\eta)=\frac{2(u-\eta)}{\eta}R_{\eta,\xi}^{23}~,
$$
Therefore
$$
R_{\eta,\xi}^{23}(\eta)L^{1}(v)({\rm qdet}_{\eta,\xi}\, L(u))
R_{\eta,\xi}^{23}(\eta)=
R_{\eta,\xi}^{23}(\eta)({\rm qdet}_{\eta,\xi}\, L(u))L^{1}(v)
R_{\eta,\xi}^{23}(\eta))~,
$$
and ${\rm qdet}_{\eta\xi}\, L(u)$ commutes with $L(v)$, i.e.
$$
[{\rm qdet}_{\eta,\xi}\, L(u),\,L(v)]=0
$$
The proof is complete.

At last we have the following theorem.
\begin{theorem}
Let $A$ be a Hopf algebra generated by matrix elements of
the generating function $L(u)$ satisfying the following conditions:
\bn
R_{\eta,\xi}(u-v)\cdot\stackrel{1}{L}(u)\stackrel{2}{L}(v)=
\stackrel{2}{L}(v)\stackrel{1}{L}(u)\cdot R_{\eta,\xi}(u-v)~,
\label{RTT11}
\ed
\bn
L(u)=\left(\begin{array}{cc}
T_{\a}^{\frac{1}{2}} & 0\\
\xi h_{\a}T_{\a}^\frac{1}{2} & T^{-\frac{1}{2}}\end{array}\right)+
\frac{L_{(0)}}{u}+\frac{L_{(1)}}{u^2}+~\cdots~,
\label{RTT12}
\ed

\bn
{\rm qdet}_{\eta,\xi}\, L(u)=1~
\label{RTT13}
\ed
\bn
\Delta(e_{ij}(u))=\sum_k e_{ik}(u)\otimes e_{kj}(u)~,
\label{RTT14}
\ed
\bn
S(L(u))=L^{-1}(u)~,
\label{RTT15}
\ed
where $T_{\a}:=1-2\xi e_{-\a}$.
Then $A$ is isomorphic to $Y_{\eta,\xi}(sl_2)$.
\label{RTTT1}
\end{theorem}
{\it Proof}. The proof is analogous to that of Proposition
\ref{RTTP1}.
Taking into account that
$$
L(u)=\Bigl(\rho^{(1)}(u)\otimes {\rm id}\Bigr)({\cal R}^{(F)})
$$
we see that there exists a homomorphism from $A$ to
$Y_{\eta,\xi}(sl_2)$. To see that this is an epimorphism, one can
use the formula (\ref{RY2}). To prove that this homomorphism is a
monomorphism one can use the same arguments as for non-deformed
Yangians (see for instance \cite{DF}).

{\it Remark.} To find the constant term in the decomposition of $L(u)$
into the series in $u^{-1}$, one should find
\bn
\lim_{u\to\infty}(\rho^{(1)}(u)\otimes {\rm id})({\cal R}^{(F)})=
(\rho^{(1)}\otimes {\rm id})({\cal R}^{(F)})=
\left(\begin{array}{cc}1 & 0\\
\xi h_{\a} & 1\end{array}\right)
\left(\begin{array}{cc}z^{-1} & 0 \\
0 & z\end{array}\right)~,
\label{RTT16}
\ed
where
\bn
z=1+\xi e_{-\a}+\frac{3!!}{2!}\xi^2 e_{-\a}^{2}+\cdots+\frac{(2n-
1)!!}{n!}
\xi^n e_{-\a}^n+\cdots~~.
\label{RTT17}
\ed
It is not difficult to see that $z^2=T_{\a}^{-1}$, hence
\bn
\lim_{u\to\infty}(\rho^{(1)}(u)\otimes {\rm id})({\cal R}^{(F)})=
\left(\begin{array}{cc}T_{\a}^\frac{1}{2} & 0\\
\xi h_{\a}T_{\a}^\frac{1}{2}&T_{\a}^{-\frac{1}{2}}\end{array}\right)~.
\label{RTT18}
\ed

\setcounter{equation}{0}
\section{Realization of the deformed Yangian $Y_{\eta,\xi}(sl_2)$
in terms of generating functions}

The realization of the usual Yangian $Y_{\eta}(sl_2)$ in terms of
the generating functions ("fields" realization) $h^{+}_{\a}(u)$ and
$e^{+}_{\pm\a}(u)$ (see Section 5) can be obtained from
the Gauss decomposition of the L-operator:
\bn
L(u)=(\rho^{(1)}(u)\otimes {\rm id})({\cal R}_{+}{\cal R}_{0}{\cal
R}_{-})
=\left(\begin{array}{cc} 1 & 0\\
-e_{\a}^+(u) & 1 \end{array}\right)
\left(\begin{array}{cc} k_1(u) & 0\\
 0 & k_2(u) \end{array}\right)
\left(\begin{array}{cc} 1 & -e_{-\a}^+(u)\\
0 & 1\end{array}\right)~.
\label{RY1}
\ed
where $k_1(u)k_2(u-1)=1$, $k_2(u)k_1^{-1}(u)=h^+(u)$
(see \cite{KT} and {\cite{DF}).

The same procedure can be applied to ${\cal R}^{(F)}$. We have
\[
L(u)=(\rho^{(1)}(u)\otimes {\rm id})(F^{21}
{\cal R}_+{\cal R}_0{\cal R}_-F^{-1})=
\]
\bn
=\left(\begin{array}{cc} 1 & 0\\
\xi h_{\a} & 1 \end{array}\right)
\left(\begin{array}{cc} 1 & 0\\
-e_{\a}^+(u) & 1\end{array}\right)
\left(\begin{array}{cc} k_1(u) & 0\\
0 & k_2(u)\end{array}\right)
\left(\begin{array}{cc} 1 & -e_{-\a}^+(u)\\
0 & 1\end{array}\right)
\left(\begin{array}{cc} T_{\a}^\frac{1}{2} & 0\\
0 & T_{\a}^{-\frac{1}{2}}\end{array}\right)=
\label{RY2}
\ed
\bn
=\left(\begin{array}{cc} 1 & 0\\
\xi h_{\a}-e_{\a}^+(u) & 1\end{array}\right)
\left(\begin{array}{cc} k_1(u)T_{\a}^\frac{1}{2} & 0\\
0 & k_2(u)T_{\a}^{-\frac{1}{2}}\end{array}\right)
\left(\begin{array}{cc} 1 &-T_{\a}^{-\frac{1}{2}}\\
0 & 1\end{array}\right)~.
\label{RY3}
\ed

\noindent
{\it Remark.} Strictly speaking, the decomposition (\ref{RY3}) is
valid in the algebra $A$ defined by the relation (\ref{RTT7}),
but simultaneously the formula (\ref{RY3}) shows that
all the generators $\{h_{\a},~e_{\a},~ T_{\a}^{\pm 1},~e_{\d-\a}\}$
of the deformed Yangian $Y_{\eta,\xi}(sl_2)$ can be expressed
in terms of the L-operator (\ref{RTT3}), what proves that the
homorphism
constructed in Theorem \ref{RTTT1} is an epimorphism.

The Gauss decomposition (\ref{RY3}) provides the following choice of
generators for $Y_{\eta,\xi}(sl_2)$:
\bn
\tilde{h}_{\a}^+(u)=
T_{\a}^{-\frac{1}{2}} h_{\a}^+(u)T_{\a}^{-\frac{1}{2}}~. \hsp
\tilde{e}_{\a}^+(u)=e_{\a}^+(u)-\xi h_{\a}~,\hsp
\tilde{e}_{-\a}^+(u)=
T_{\a}^{-\frac{1}{2}} e_{-\a}^+(u)T_{\a}^{-\frac{1}{2}}~.
\label{RY4}
\ed
Using relations (\ref{PS4})-(\ref{PS10}) one can obtain
the following relations between
$\tilde{h}_{\a}(u),~\tilde{e}_{\pm\a}(u)$ and
$T_{\a}^{\pm\frac{1}{2}}$:
\bn
T_{\a}^{\pm\frac{1}{2}}\tilde{h}_{\a}^+(u)T_{\a}^{\mp\frac{1}{2}}=
\Bigl(1\pm\eta\xi\tilde{e}_{-\a}^+(u)\Bigr)^{\mp1}
\tilde{h}_{\a}^+(u)\Bigr(1\pm\eta\xi\tilde{e}_{-\a}^+(u)\Bigr)^{-1},
\label{RY5}
\ed
\bn
T_{\a}^{\pm 1}\tilde{e}_{\a}^+(u)T_{\a}^{\mp1}=
\tilde{e}_{\a}^+(u)\mp2\eta\xi\pm
2\eta\xi\Bigl(1\pm\eta\xi\tilde{e}_{-\a}^+(u)\Bigr)^{\pm 1}
\tilde{h}_{\a}^+(u)\Bigl(1\pm\eta\xi\tilde{e}_{-\a}^+(u)\Bigr)^{-1},
\label{RY6}
\ed
\bn
T_{\a}^{\pm\frac{1}{2}}\tilde{e}_{-\a}(u)T_{\a}^{\mp\frac{1}{2}}=
\tilde{e}_{-\a}(u)\Bigl(1\pm\eta\xi\tilde{e}_{-\a}(u)\Bigr)^{-1}.
\label{RY7}
\ed
With the help of the relations (\ref{PS4})-(\ref{PS10}) and
(\ref{RY5})-(\ref{RY7}) we can prove the following theorem.
\begin{theorem}
Defining relations for the "fields" $\tilde{h}_{\a}^+(u)$,
$\tilde{e}_{\pm\a}^+(u)$ have the form:
\bn
{\cal H}_1(u)\tilde{h}_{\a}^+(v)=
{\cal H}_1(v)\tilde{h}_{\a}^+(u)~,
\label{RY8}
\ed
\bn
[\tilde{e}_{\a}^+(u),\tilde{e}_{\a}^+(v)]=
\frac{\eta\Bigl(\tilde{e}_{\a}^+(v)-\tilde{e}_{\a}^+(u)\Bigr)^2}{v-u}+
2\eta\xi\Bigr(\tilde{e}_{\a}^+(u)-\tilde{e}_{\a}^+(v)\Bigl)~,
\label{RY9}
\ed
\bn
(u-v+\eta){\cal G}_1(u)\tilde{e}_{-\a}^+(v)-
(u-v-\eta){\cal G}_1(v)\tilde{e}_{-\a}^+(u)=
\eta \Bigl({\cal G}_1(u)\tilde{e}_{-\a}^+(u))+
{\cal G}_1(v)\tilde{e}_{-\a}^+(v)\Bigr)~,
\label{RY10}
\ed
\bn
(u-v-\eta){\cal H}_1(u)\tilde{e}_{-\a}^+(v)-
(u-v+\eta){\cal G}_1(v)\tilde{h}_{\a}^+(u)=
\eta {\cal H}_1(v)\tilde{e}_{-\a}^+(u)+
\eta {\cal G}_1(u)\tilde{h}_{\a}(v)~,
\label{RY11}
\ed
\bn
\Bigl(\tilde{e}_{\a}^+(u)-2\eta\xi+
2\eta\xi {\cal H}_2(u)\Bigl){\cal G}_2(v)-
{\cal G}_2(v)\tilde{e}_{\a}^+(u)
=-\frac{\eta\Bigl({\cal H}_2(u)-
{\cal H}_2(v)\Bigr)}{u-v}+
2\eta\xi {\cal G}_2(v)~,
\label{RY12}
\ed
\[
(u-v+\eta)\Bigl(\tilde{e}_{\a}^+(v)-2\eta\xi+
2\eta\xi {\cal H}_2(v)\Bigr){\cal H}_2(u)-(u-v-\eta){\cal H}_2(u)=
\]
\bn
=\eta {\cal H}_2(v)\tilde{e}_{\a}^+(u)+
\eta\Bigl(\tilde{e}_{\a}^+(u)-2\eta\xi+
2\eta\xi {\cal H}_2(u)\Bigr){\cal H}_2(u)~.
\label{RY13}
\ed

where
\bn
{\cal H}_1(u)=\Bigl(1-2\eta\xi\tilde{e}_{-\a}^+(u)\Bigr)^{-1}
\tilde{h}_{\a}^+(u)\Bigl(1-2\eta\xi\tilde{e}_{-\a}^+(u)\Bigr)^{-1}~,
\label{RY14}
\ed
\bn
{\cal H}_2(u)=\Bigl(1-\eta\xi\tilde{e}_{-\a}^+(u)\Bigr)^{-1}
\tilde{h}_{\a}^+(u)\Bigl(1-\eta\xi\tilde{e}_{-\a}^+(u)\Bigr)^{-1}~,
\label{RY15}
\ed
\bn
{\cal G}_1(u)=\Bigl(1-2\eta\xi\tilde{e}_{-\a}^+(u)\Bigr)^{-1}
\tilde{e}_{-\a}^+(u)~, \hsp
{\cal G}_2(u)=\Bigl(1-\eta\xi\tilde{e}_{-\a}^+(u)\Bigr)^{-1}
\tilde{e}_{-\a}^+(u)~.
\label{RY16}
\ed
\end{theorem}
{\it Proof.} Let us prove for instance the formula {\ref{RY10}).
 From (\ref{PS9}) we have
$$
(u-v+\eta)e_{-\a}^+(u)e_{-\a}^+(v)-(u-v-\eta)e_{-\a}^+(v)e_{-\a}^+(u)=
\eta(e_{-\a}^+(u))^2+\eta(e_{-\a}^+(v))^2~.
$$
Substituting $e_{-\a}(u)=
T_{\a}^{\frac{1}{2}}\tilde{e}^{-\a}_+(u)T_{\a}^{\frac{1}{2}}$ and
using
that $T_{\a}^{-\frac{1}{2}}\tilde{e}_{-\a}^+(u)T_{\a}^{\frac{1}{2}}=
\tilde{e}_{-\a}^+(u)\Bigl(1-\eta\xi\tilde{e}_{-\a}^+(u)\Bigr)^{-1}$
(see (\ref{RY7}), we obtain the formula (\ref{RY10}).
\vspace{0.5cm}

{\bf Acknowledgments}
The authors are thankful to the Swedish Academy of Science for
the support of the visit of the first and third authors to the Royal
Institute of Technology (Stockholm), during which visit the
paper was completed. The first and third authors would like also to
thank
the Russian Foundation for Fundamental Research, grant No. 95-01-
00814A,
and the ISF, grant MBI300, for the financial support.

\end{document}